\documentclass[4pt, oneside]{article}
\usepackage{yfonts}
\usepackage{amsmath,amssymb,amsthm,units,stmaryrd,relsize}
\usepackage{graphicx}
\usepackage[all]{xy}

\newtheorem{definition}{Definition}[section]
\newtheorem{lemm}{Lemma}[section]
\newtheorem{prop}{Proposition}[section]
\newtheorem{cor}{Corollary}[section]
\newtheorem{theorem}{Theorem}[section]

\def\lb{\left\llbracket}
\def\rb{\right\rrbracket}
\def\nx{f}

\def\fmodels{\xymatrix{
\ar@{|=}[r]^{<\omega}&
}
}
\def\nmodels{\xymatrix{
\ar@{|=}[r]^{N}&
}
}
\def\formtwo{{\tt Cycle}}
\def\formthree{{\tt Start}}
\def\formfour{{\tt Bundle}}
\def\formfive{{\tt Tangle}}
\def\formsix{{\tt Trouble}}
\def\<{\left <}
\def\next{{f}}
\def\hf{{[f]}}
\def\ev{{\left< f\right>}}
\def\nc{{\Box}}
\def\ps{{\Diamond}}
\def\Ps{\mathop{\vcenter{\hbox{$\mathlarger{\mathlarger{\mathlarger{\mathlarger{\ps}}}}$}}}}

\def\peq{\preccurlyeq}

\def\>{\right >}

\def\cbra{\left \{}
\def\cket{\right \}}

\DeclareSymbolFont{AMSb}{U}{msb}{m}{n}
\DeclareMathSymbol{\N}{\mathbin}{AMSb}{"4E}
\DeclareMathSymbol{\Z}{\mathbin}{AMSb}{"5A}
\DeclareMathSymbol{\R}{\mathbin}{AMSb}{"52}
\DeclareMathSymbol{\Q}{\mathbin}{AMSb}{"51}
\DeclareMathSymbol{\I}{\mathbin}{AMSb}{"49}
\DeclareMathSymbol{\C}{\mathbin}{AMSb}{"43}
\begin{document}
\date{\today}
\title{Non-finite axiomatizability of Dynamic Topological Logic
      }
\author{David Fern\'{a}ndez-Duque
\\ {\small Sevilla}$  $
        }
\maketitle

  \begin{abstract}
  {\em Dynamic topological logic} ($\mathcal{DTL}$) is a polymodal logic designed for reasoning about {\em dynamic topological systems}. These are pairs $\langle X,f\rangle$, where $X$ is a topological space and $f:X\to X$ is continuous. $\mathcal{DTL}$ uses a language $\mathsf L$ which combines the topological $\mathsf{S4}$ modality $\nc$ with temporal operators from linear temporal logic.
  
  Recently, I gave a sound and complete axiomatization $\mathsf{DTL}^\ast$ for an extension of the logic to the language $\mathsf L^\ast$, where $\ps$ is allowed to act on finite sets of formulas and is interpreted as a tangled closure operator. No complete axiomatization is known over $\mathsf L$, although one proof system, which we shall call $\mathsf{KM}$, was conjectured to be complete by Kremer and Mints.
  
In this paper we show that, given any language $\mathsf L'$ such that ${\sf L}\subseteq {\sf L}'\subseteq{\sf L}^\ast$, the set of valid formulas of ${\sf L}'$ is not finitely axiomatizable. It follows, in particular, that ${\sf KM}$ is incomplete.
  \end{abstract}

\section{Introduction}

Finding a transparent axiomatization for Dynamic Topological Logic ($\mathcal {DTL}$) has been an elusive open problem since 2005, when one (which we shall call $\mathsf{KM}$) was proposed by Kremer and Mints in \cite{kmints} without establishing its completeness. In \cite{dtlaxiom} I offered a complete axiomatization, not over the language $\mathsf L$ used in \cite{kmints}, but rather in an extended language $\mathsf L^\ast$ which allowed the modal $\ps$ to be applied to finite sets of formulas. It was then interpreted as a `tangled closure' operator (see Section \ref{sectan}). The resulting logic is called $\mathsf{DTL}^\ast$.

However, the fact that $\mathsf{DTL}^\ast$ used the unfamiliar `tangled closure' operation and was substantially less intuitive than $\mathsf{KM}$ left the completeness of the latter as a relevant open problem. Actually, the only motivation given in \cite{dtlaxiom} for passing to an extended language was, {\em There is a completeness proof which works in the extended language but not in the original one}; a valid, but not terribly compelling, reason.

The results in this paper will show that indeed the use of the tangled closure is an essential part of this axiomatization, and cannot be removed without extending $\mathsf{KM}$, although it is not clear what such an extension should look like. In fact, we prove more. We show that, given $k<\omega$, there is a formula $\formsix^k\in \mathsf L$ such that $\formsix^k$ is derivable in $\mathsf{DTL}^\ast$ only by using formulas of the form $\ps\Gamma$, where $\Gamma$ has at least $k$ elements. This shows that $\cal DTL$ can be written as a strictly increasing sequence of theories and hence is not finitely axiomatizable; it follows, in particular, that $\mathsf{KM}$ is incomplete.

\subsection{Previous work on $\mathcal{DTL}$}

{\em Dynamic topological logic} ($\mathcal {DTL}$) combines the topological $\mathsf{S4}$ with Linear Temporal Logic. The `topological interior' interpretation of modal logic was already studied by Tarski, McKinsey and others around the 1940s \cite{tarski} and more recently in works like \cite{s4real,kstrong,mints}. Temporal logic also has a long history, having been studied by Prior before 1960 \cite{prior} and received substantial attention since; see \cite{temporal} for a nice overview.

The purpose of $\mathcal{DTL}$ is to reason about {\em dynamic topological systems} (dts's); these are pairs $\langle X,f\rangle$, where $X$ is a topological space and $f:X\to X$ is a continuous function, and shall be discussed in greater detail in Section \ref{basic}. Dynamic Topological Logic was originally introduced as a bimodal logic in \cite{arte}, where it was called $\mathsf{S4C}$. In our notation, it uses the `interior' modality $\nc$, interpreted topologically, and `next-time' modality $\nx$, interpreted as a preimage operator; see Section \ref{basic} for details. The logic $\mathsf{S4C}$ is a rather well-behaved modal logic; it is decidable, axiomatizable and has the finite model property, all of which was established in \cite{arte}. Later, \cite{kmints} showed that a variant, called $\mathsf{S4H}$, is complete for the class of dynamical systems where $f$ is a homeomorphism.

Also in \cite{kmints}, it was noted that by adding the infinitary temporal modality `henceforth' (here denoted $\hf$), one could reason about long-term behavior of dts's, capturing phenomena such as topological recurrence. Thus they introduced an extension of $\mathsf{S4C}$, which we denote $\mathcal{DTL}$. $\mathcal{DTL}$ turned out to behave much worse than $\mathsf{S4C}$; it was proven to be undecidable in \cite{konev}, and in \cite{wolter} it was also shown that, if we restrict to the case where $f$ is a homeomorphism, then the logic becomes non-axiomatizable. Fortunately, with arbitrary continuous functions the logic turned out to be recursively enumerable \cite{me2}, but the only currently known axiomatization is over $\mathsf L^\ast$ \cite{dtlaxiom}.

This axiomatization uses the fact, first observed in \cite{do}, that $\mathsf L^\ast$ is more expressive than $\mathsf{L}$. There, it is shown that, over the class of finite $\mathsf{S4}$ models, $\mathsf L^\ast_\ps$ (i.e., the fragment of $\mathsf L^\ast$ without temporal modalities) is equally expressive to the bisimulation-invariant fragment of both first-order logic and monadic second-order logic, while $\mathsf L_\ps$ is weaker. This added expressive power is used in an important way in \cite{dtlaxiom}, although it is not proven that such an extension is {necessary}.

\subsection{Layout of the paper}

Sections \ref{topre}-\ref{theax} give a general review of dynamic topological logic and the known results relevant for this paper. In Section \ref{topre}, we introduce topological spaces and show how one can see preorders as a special case. Section \ref{sectan} gives the main properties of the tangled closure operator, an important part of $\mathcal{DTL}^\ast$, introduced in Section \ref{basic}. Section \ref{theax} then reviews the axiomatization from \cite{dtlaxiom} and defines some important sublogics.

Section \ref{tb} introduces tangled bisimulations, which are based on those presented in \cite{do}. These are used to show that $\Ps_{i=1}^{k+1}\gamma_i$ cannot in general be defined using exclusively $k$-adic occurrences of $\ps$.

Section \ref{formulas} defines the formulas $\formsix^k$ which are derivable in $\mathsf{DTL}^k$, as well as other formulas which are useful for our purposes. Finally, Section \ref{secinc} shows that, indeed, $\mathsf{DTL}^k$ is consistent with $\neg\formsix^{k+1}$, thus stratifying Dynamic Topological Logic into a strictly increasing sequence of theories, from which it follows that it is not finitely axiomatizable.

\section{Topologies and preorders}\label{topre}

In this section we shall very briefly review some basic notions from topology. As is well-known, topological spaces provide an interpretation of the modal logic $\mathsf{S4}$, generalizing its familiar Kripke semantics.

Let us recall the definition of a {\em topological space}:

\begin{definition}
A {\em topological space} is a pair $\mathfrak X=\<|\mathfrak X|,\mathcal{T}_\mathfrak X\>,$ where $|\mathfrak X|$ is a set and $\mathcal T_\mathfrak X$ a family of subsets of $|\mathfrak X|$ satisfying
\begin{enumerate}
\item $\varnothing,|\mathfrak X|\in \mathcal T_\mathfrak X$;
\item if $U,V\in \mathcal T_\mathfrak X$ then $U\cap V\in \mathcal T_\mathfrak X$ and
\item if $\mathcal O\subseteq\mathcal T_\mathfrak X$ then $\bigcup\mathcal O\in\mathcal T_\mathfrak X$.
\end{enumerate}

The elements of $\mathcal T_\mathfrak X$ are called {\em open sets}. Complements of open sets are {\em closed sets}.
\end{definition}

Given a set $A\subseteq |\mathfrak X|$, its {\em interior}, denoted $A^\circ$, is defined by
\[A^\circ=\bigcup\cbra U\in\mathcal T_\mathfrak X:U\subseteq A\cket.\]
Dually, we define the closure $\overline A$ as $|\mathfrak X|\setminus(|\mathfrak X|\setminus A)^\circ$; this is the smallest closed set containing $A$.

Topological spaces generalize transitive, reflexive Kripke frames. Recall that these are pairs $\mathfrak W=\<|\mathfrak W|,\peq_\mathfrak W\>$, where $\peq_\mathfrak W$ is a preorder on the set $|\mathfrak W|$. We will write $\peq$ instead of $\peq_\mathfrak W$ whenever this does not lead to confusion.

To see a preorder as a special case of a topological space, define
\[\mathop\downarrow w=\cbra v:v\peq w\cket.\]
Then consider the topology $\mathcal T_\peq$ on $|\mathfrak W|$ given by setting $U\subseteq|\mathfrak W|$ to be open if and only if, whenever $w\in U$, we have $\mathop\downarrow w\subseteq U$ (so that the sets of the form $\downarrow w$ provide a basis for $\mathcal T_\peq$). A topology of this form is a {\em preorder topology}\footnote{Or, more specifically, a {\em downset topology}. Note that I stray from convention, since most authors use the upset topology here, but I find this presentation more natural. This will later be reflected in the semantics for $\nc$.}.

Throughout this text we will often identify preorders with their corresponding topologies, and many times do so tacitly.

We will also use the notation
\begin{itemize}
\item $w\prec v$ for $w\peq v$ but $v\not\peq w$ and
\item $w\approx v$ for $w\peq v$ and $v\peq w$.
\end{itemize}

The relation $\approx$ is an equivalence relation; the equivalence class of a point $x\in |\mathfrak W|$ is usually called a {\em cluster}, and we will denote it by $[x]$.

\section{The tangled closure}\label{sectan}

The {\em tangled closure} is an important component of $\mathsf{DTL}^\ast$. It was introduced in \cite{do} for Kripke frames and has also appeared in \cite{me:simulability,me:tangle,dtlaxiom,dynamictangle}.

\begin{definition}
Let $\mathfrak X$ be a topological space and $\mathcal S\subseteq 2^{|\mathfrak X|}$.

Given $E\subseteq|\mathfrak X|$, we say $\mathcal S$ is {\em tangled} in $E$ if, for all $A\in \mathcal S$, $A\cap E$ is dense in $E$.

We define the {\em tangled closure of $\mathcal S$}, denoted ${\mathcal S}^\ast$, to be the union of all sets $E$ such that $\mathcal S$ is tangled in $E$.
\end{definition}

It is important for us to note that the tangled closure is defined over any topological space; however, we will often be concerned with locally finite preorders in this paper. Here, the tangled closure is relatively simple.

\begin{lemm}
Let $\<S,\peq\>$ be a finite preorder, $x\in S$ and $\mathcal O\subseteq 2^S$. Then, $x\in{\mathcal O}^\ast$ if and only if there exist $\<y_A\>_{A\in\mathcal O}$ such that $y_A\in A$, $y_A\peq x$ for all $A\in\mathcal O$ and $y_A\approx y_{B}$ for all $A,B\in \mathcal O$.
\end{lemm}

\proof
A proof can be found in any of \cite{me:simulability,me:tangle,dtlaxiom,dynamictangle}.
\endproof

\section{Dynamic Topological Logic}\label{basic}

The language $\mathsf L^\ast$ is built from propositional variables in a countably infinite set $\mathsf{PV}$ using the Boolean connectives $\wedge$ and $\neg$ (all other connectives are to be defined in terms of these), the unary modal operators $\next$ (`next') and $\hf$ (`henceforth'), along with a polyadic modality $\ps$ which acts on finite sets, so that if $\Gamma$ is a finite set of formulas then $\ps\Gamma$ is also a formula. Note that this is a modification of the usual language of $\mathcal{DTL}$, where $\ps$ acts on single formulas only. We write $\nc$ as a shorthand for $\neg\ps\neg$; similarly, $\ev$ denotes the dual of $\hf$. We also write $\ps\gamma$ instead of $\ps\cbra\gamma\cket$; its meaning is identical to that of the usual $\mathsf{S4}$ modality \cite{dynamictangle}. We will often write $\Ps_{n =1}^N\gamma_i$ instead of $\ps\{\gamma_i\}_{1\leq i\leq N}$. 

Given a formula $\phi$, the {\em depth} of $\phi$, denoted $\mathrm{dpt}(\phi)$, is the modal nesting depth of $\phi$, while its {\em width}, $\mathrm{wdt}(\phi)$, denotes the maximal $k$ such that $\phi$ has a subformula of the form $\Ps_{i=1}^k\gamma_i$. For $k<\omega$, $\mathsf L^k$ denotes the sublanguage of $\mathsf L^\ast$ where all formulas have width at most $k$. Thus, in particular, $\mathsf L=\mathsf L^1$.

Formulas of $\mathsf L^\ast$ are interpreted on dynamical systems over topological spaces, or {\em dynamic topological systems}.

\begin{definition}
A {\em weak dynamic topological system (dts)} is a triple
\[\mathfrak X=\<|\mathfrak X|,\mathcal{T}_{\mathfrak X},f_{\mathfrak X}\>,\]
where $\<|\mathfrak X|,\mathcal{T}_{\mathfrak X}\>$ is a topological space and
\[f_{\mathfrak X}:|\mathfrak X|\to |\mathfrak X|.\]

If further $f_\mathfrak X$ is continuous\footnote{That is, whenever $U\subseteq|\mathfrak X|$ is open, then so is $f^{-1}(U)$}, we say $\mathfrak X$ is a {\em dynamical system}.
\end{definition}

\begin{definition}
Given a (weak) dynamic topological system $\mathfrak X$, a {\em valuation} on $\mathfrak X$ is a function
\[\lb\cdot\rb:\mathsf L^\ast\to 2^{|\mathfrak X|}\]
satisfying
\begin{align*}
\lb\alpha\wedge\beta\rb_{\mathfrak X}&=\lb\alpha\rb_{\mathfrak X}\cap \lb\beta\rb_{\mathfrak X}\\
\lb\neg\alpha\rb_{\mathfrak X}&=|\mathfrak X|\setminus \lb\alpha\rb_{\mathfrak X}\\
\lb\next\alpha\rb_{\mathfrak X}&=f^{-1}\lb\alpha\rb_{\mathfrak X}\\
\lb\hf\alpha\rb_{\mathfrak X}&=\displaystyle\bigcap_{n\geq 0}f^{-n}\lb\alpha\rb_{\mathfrak X}\\
\lb\ps\cbra\alpha_1,...,\alpha_n\cket\rb_\mathfrak X&=\cbra\lb\alpha_1\rb_\mathfrak X,...,\lb\alpha_n\rb_\mathfrak X\cket^\ast.
\end{align*}
\end{definition}

A {\em (weak) dynamic topological model} (wdtm/dtm) is a (weak) dynamic topological system $\mathfrak X$ equipped with a valuation $\lb\cdot\rb_\mathfrak X$. We say a formula $\phi$ is {\em valid} on $\mathfrak X$ if $\lb\phi\rb_\mathfrak X=|\mathfrak X|$, and write $\mathfrak X\models\phi$. If a formula $\phi$ is valid on every dynamic topological model, then we write $\models\phi$. $\mathcal{DTL}$ is the set of valid formulas of $\sf L$ under this interpretation, while $\mathcal{DTL}^\ast$ denotes the set of valid formulas of ${\sf L}^\ast$.

We will often write $\<\mathfrak X,x\>\models\phi$ instead of $x\in\lb\phi\rb_{\mathfrak X}$.

\section{The axiomatization}\label{theax}

We shall distinguish $\mathcal{DTL}^\ast$, which is defined semantically, from $\mathsf{DTL}^\ast$, which is a proof system. The two have the same set of theorems, but we will be interested in natural subsystems of $\mathsf{DTL}^\ast$ which are defined syntactically.

Below, note that the modality $f$ is unary, and $f\Gamma$ is merely a shorthand for $\{f\gamma:\gamma\in\Gamma\}$; $p$ denotes a propositional variable and $P$ a finite set of propositional variables. Then, the axiomatization $\mathsf{DTL}^\ast$ consists of the following:

\begin{description}
\item[$\mathsf{Taut}$] All propositional tautologies.
\item[Topological axioms]\
\begin{description}
\item[$\mathsf K$] $\nc(p\to q)\to(\nc p\to\nc q)$
\item[$\mathsf{T}$] $\nc p\to p$
\item[$\mathsf 4$] $\nc p\to \nc\nc p$
\item[$\mathsf{Fix}_\ps$] $\ps P\to\bigwedge_{q\in P}\ps(q\wedge\ps P)$
\item[$\mathsf{Ind}_\ps$] ${\nc\displaystyle\bigwedge_{q\in P} (p\to\ps(q\wedge p))}\to
{(p\to\displaystyle\ps P)}$
\end{description}
\item[Temporal axioms]\
\begin{description}
\item[$\mathsf{Neg}_\next$]$\neg\next p\leftrightarrow\next\neg p$
\item[$\mathsf{And}_\next$] $\next(p\wedge q)\leftrightarrow \next p\wedge \next q$
\item[$\mathsf{Fix}_\hf$] $\hf p\to p\wedge\next\hf p$
\item[$\mathsf{Ind}_\hf$]
${\hf(p\to\nx p)}\to({p\to\hf p})$
\end{description}
\item[$\mathsf{Cont}^\ast$]$\ps\next P\to\next\ps P$.
\item[Rules]\
\begin{description}
\item[$\mathsf{MP}$] Modus ponens
\item[$\mathsf{Subs}$] $\displaystyle\dfrac{\phi}{\phi[\vec p/\vec\psi]}$
\item[$\mathsf N_\nc$] $\displaystyle\dfrac{\phi}{\nc\phi}$\hskip 60pt$\mathsf N_\next$ $\displaystyle\dfrac{\phi}{\nx\phi}$\hskip 60pt$\mathsf N_\hf$  $\displaystyle\dfrac{\phi}{\hf\phi}$
\end{description}
\end{description}

This axiomatization is sound and complete, as proven in \cite{dtlaxiom}:

\begin{theorem}
$\mathsf{DTL}^\ast$ is sound and complete for the class of dynamic topological models.
\end{theorem}

There are many subtleties in our proof system, so before continuing we should make a few remarks.

First, let us say a few words about the substitution rule. It is to be understood as `simulataneous substitution', where $\vec p$ represents a finite sequence of variables, $\vec\psi$ a finite sequence of formulas and each variable is replaced by the respective formula. By standard arguments, this rule preserves validity, as there is nothing in our semantics distinguishing atomic facts from complex propositions.

Further, since we are concerned with finite axiomatizability of a logic it is important to include it; otherwise, each substitution instance of any of the axioms would have to be regarded as a new axiom and the finite axiomatizability would fail for obvious reasons. Of course this is not the only possible presentation, as one can also consider axiomatizations by finitely many {\em schemas}, but here we shall consider different formulas to be different also as axioms.

With this in mind, we should note that the above axiomatization is not finite, nor can it be modified into a finite version in an obvious way. Evidently the set of all propositional tautologies can be replaced by finitely many axioms, but this is not what concerns us. Much more importantly, we need infinitely many axioms for $\ps$, and it is only in the metalanguage that we can give them a uniform presentation. In fact, the symbol $P$ representing a finite set of propositional variables is not a symbol of $\sf L^\ast$, where we would have to write out explicitly $\{p_1,\hdots, p_k\}$ for each given value of $k$.

Of particular interest is the schema $\mathsf{Cont}^\ast$. This was originally named $\mathsf{TCont}$; we adopt the new notation to stress that the standard `continuity' axiom,
\[\mathsf{Cont}^1:\ps\nx p\to\nx\ps p,\]
is indeed a special case.

${\sf Cont}^\ast$ is really an infinite collection of axioms. To be precise, for $k<\omega$ let
\[\mathsf{Cont}^k=\Ps_{i\in [1,k]}\nx p_i\to \nx\Ps_{i\in [1,k]} p_i.\]
Note that $\mathsf{Cont}^{k+1}$ extends $\mathsf{Cont}^{k}$ since we can always substitute $p_{k+1}$ by $p_k$. 

We then let $\mathsf{DTL}^k$ be the variant of $\mathsf{DTL}^\ast$ where $\mathsf{Cont}^\ast$ is replaced by $\mathsf{Cont}^k$. We denote derivability in $\mathsf{DTL}^k$ by $\vdash^k$. $\mathsf{DTL}^0$ denotes the system with no continuity axiom.

Our goal will be to show that $\langle\mathsf{DTL}^k\rangle_{k<\omega}$ gives a sequence of theories of strictly increasing strength. Since ${\sf DTL}^\ast$ is the union of these theories, it will follow as a straightforward consequence that ${\sf DTL}^\ast$ is not finitely axiomatizable. However, to do this we will need a second refinement, this time of each ${\sf DTL}^k$.

For $n,k<\omega$, we let ${\sf DTL}_n^{k}$ be the subtheory of ${\sf DTL}^k$ which restricts the substitution rule in the following ways:
\begin{enumerate}
\item ${\sf Subs}$ may only be applied immediately to axioms and
\item if $\sf Subs$ is applied to ${\sf Cont}^k$, then each $p_i$ must be replaced by a formula with modal depth at most $n$.
\end{enumerate}

A very easy induction on derivations shows that any proof in ${\sf DTL}^k$ may be transformed into one satisfying the above two conditions for some value of $n$ and hence ${\sf DTL}^k=\bigcup_{n<\omega}{\sf DTL}^k_n$. We denote derivability in ${\sf DTL}^k_n$ by $\vdash^k_n$.

The reason for passing to ${\sf DTL}^k_n$ is that the substitution rule, while preserving validity, does not preserve {\em model} validity; if $\mathfrak M\models\phi$, it does not always follow that $\mathfrak M\models \phi[p/\psi]$. Later we wish to build specific models of fragments of ${\sf DTL}^\ast$, and to check soundness for these models, ${\sf DTL}^k_n$ has the advantage that we only need to focus on substitution instances of axioms. This will become relevant in Section \ref{secinc}.

${\sf DTL}^\ast$ is an extension of ${\sf KM}$, which can be defined as follows:
\begin{definition}
The calculus $\mathsf{KM}$ is the restriction\footnote{Of course this description is anachronical, and it would be more accurate to think of $\mathsf{DTL}^\ast$ as an extension of $\sf KM$.} of $\mathsf{DTL}^\ast$ to $\mathsf{L}^1$.
\end{definition}
 In $\sf KM$, all appearances of $\ps$ must be applied to a single formula; in particular, the axioms ${\sf Fix}_\ps$ and ${\sf Ind}_\ps$ are not present, and ${\sf Cont}^\ast$ becomes ${\sf Cont}^1$. We should note that $\mathsf{DTL}^1$ is very similar, but not identical, to $\mathsf{KM}$.  $\mathsf{DTL}^1$ allows formulas of the form $\ps\Gamma$ within derivations for $\Gamma$ arbitrarily large, but $\mathsf{Cont}^\ast$ is also replaced by $\mathsf{Cont}^1$. We do have, however, that ${\sf KM}\subseteq{\sf DTL}^1$.

Later we shall show that the sequence $\langle \mathsf{DTL}^k\rangle_{k<\omega}$ is strictly increasing in strength, even over $\mathsf L$; i.e., there are formulas $\formsix^{k}\in{\sf L}$ such that $\vdash^{k+1}\formsix^{k+1}$ but $\not\vdash^k\formsix^{k+1}$. These are defined in Section \ref{formulas}; but first, we need to define {\em partial tangled bisimulations}, the fundamental tool we shall use to prove our main results.

\section{Tangled bisimulations}\label{tb}

Our main results are based on partial bisimulation techniques. As we will be working in a polyadic system, we shall need a notion of partial bisimulation which preserves the polyadic $\ps$. Such a notion was already introduced in \cite{do}; here we present a slight generalization which is more sensitive to the width of formulas. For more information on partial bisimulations, we refer the reader to a text such as \cite{black}.

\begin{definition}[Tangled partial bisimulation]

Given models $\mathfrak X,\mathfrak Y$, $n<\omega$ and $k\leq\omega$, we define a binary relation $\leftrightarroweq^n_k\subseteq|\mathfrak X|\times|\mathfrak Y|$ by inducion on $n$ as follows.

For $n=0$, $x\leftrightarroweq^0_k y$ if and only if $x$ and $y$ satisfy the same set of atoms.

Otherwise, $x\leftrightarroweq^{n+1}_k y$ if $x,y$ satisfy the same set of atoms and
\begin{description}
\item[$\mathsf{Forth}_\peq$] whenever $m< k$ and $x_1\approx x_2\approx\hdots\approx x_{m}\peq x$, there are $y_1\approx y_2\approx\hdots y_{m}\peq y$ such that $x_i\leftrightarroweq^n_k y_i$ for all $i\leq m$,
\item[$\mathsf{Back}_\peq$] whenever $m<k$ and $y_1\approx y_2\approx\hdots\approx y_{m}\peq y$, there are $x_1\approx x_2\approx\hdots \approx x_{m}\peq x$ such that $x_i\leftrightarroweq^n_k y_i$ for all $i\leq m$,
\item[$\mathsf{Forth}_\nx$] $f_\mathfrak X(x)\leftrightarroweq^n_k f_\mathfrak Y(y)$
\item[$\mathsf{Forth}_\hf$] for every $m<\omega$ there is $m'<\omega$ such that $f^m_\mathfrak X(x)\leftrightarroweq^n_k f^{m'}_\mathfrak Y(y)$ and
\item[$\mathsf{Back}_\hf$] for every $m<\omega$ there is $m'<\omega$ such that $f^m_\mathfrak Y(y)\leftrightarroweq^n_k f^{m'}_\mathfrak X(x)$.
\end{description}
\end{definition}

We will write $\leftrightarroweq^n_\ast$ when $k=\omega$; in this case there are no bounds on the clauses for $\peq$. Note that there is no `back' clause for $\nx$ as it would be identical to $\mathsf{Forth}_\nx$. For purely topological structures (i.e., without the function $f_\mathfrak X$), we shall also use the analogous notion of partial bisimulation, simply removing the clauses for $\nx,\hf$.

When the respective structures are clear from context, we may write $x\leftrightarroweq^n_k y$ instead of $\langle\mathfrak X,x\rangle\leftrightarroweq^n_k\langle\mathfrak Y,y\rangle.$

\begin{lemm}\label{basicbis}
If $\varphi$ is a formula with $\mathrm{dpt}(\varphi)\leq n$ and $\mathrm{wdt}(\varphi)< k$ and $\mathfrak X,\mathfrak Y$ are finite dtm's, then whenever $\langle\mathfrak X,x\rangle\leftrightarroweq^{n}_k \langle\mathfrak Y,y\rangle$, we have that $x\in\lb\varphi\rb_\mathfrak X$ if and only if $y\in\lb\varphi\rb_\mathfrak Y$.
\end{lemm}

\proof
The proof proceeds by a standard induction on $\mathrm{dpt}(\varphi)$ and we omit it.
\endproof

Below and throughout the text, $\coprod$ denotes a disjoint union.

\begin{definition}[Simple models]
Let $k<\omega$.

A preordered model $\mathfrak S$ is {\em $k$-simple} (or merely {\em simple}) if
\[|\mathfrak S|=\coprod_{i=1}^k\lb p_i\rb_\mathfrak S.\]

If $\mathfrak S$ is $k$-simple and $x\in|\mathfrak S|$, we write $p_\mathfrak S(x)$ for the unique $p\in\{p_1,\hdots,p_k\}$ such that $x\in\lb p\rb_\mathfrak S$.
\end{definition}

As always, we will drop subindices when it does not lead to confusion, writing $p(x)$ instead of $p_\mathfrak S(x)$.

Before continuing, let us establish a notational convention. Given natural numbers $n,k$, we will denote by $|n|_k$ the unique element $m$ of $\{1,\hdots,k\}$ such that $n\equiv m\pmod{k}$. Note that this strays from the standard remainder in that $|k|_k=k$, but it shall simplify several expressions later on. Intervals shall be assumed to be intervals of natural numbers, i.e.
\[[a,b]=\{n\in \mathbb N:a\leq n\leq b\}.\]
Further, it will be convenient to assume that the set of propositional variables is enumerated by $\langle p_k\rangle_{k<\omega}$.

In the remainder of this section, we shall use partial tangled bisimulations to show that $\mathsf L^{k+1}$ is more expressive than $\mathsf L^k$. This might not be too surprising given results in \cite{do,me:simulability}, but to the best of my knowledge this has not been stated explicitly before and will provide a good warm-up for the techniques we shall use later on.

To be precise, by {\em more expressive} we mean the following:  given languages $\lambda,\lambda'\subseteq\sf L^\ast$ and a class of models $\mathcal X$, we say $\lambda'$ is {\em at least as expressive as} $\lambda$ over $\mathcal X$ if, given $\varphi\in\lambda$, there is $\varphi'\in\lambda'$ such that, for every $\mathfrak X\in\mathcal X$, $\lb\varphi\rb_\mathfrak X=\lb\varphi'\rb_\mathfrak X$. The language $\lambda'$ is {\em more expressive} than $\lambda$ if $\lambda'$ is at least as expressive as $\lambda$, but not vice-versa.

The following structures will be useful in proving our expressiveness result:

\begin{definition}[$\mathfrak A(N,K)$]\label{defank}
Given natural numbers $N,K$ we define a $K$-simple structure $\mathfrak A=\mathfrak A(N,K)$ as follows:
\begin{enumerate}
\item $|\mathfrak A|$ is the set of all pairs $(h,k)$ such that $h\in[0,NK]$, $k\in[1,K]$ and either
\begin{enumerate}
\item $h=0$ or
\item $k\not= |h|_{K}$;
\end{enumerate}
\item $(h,k)\peq_\mathfrak A(h',k')$ if and only if $h\geq h'$;
\item $p(h,k)=p_k$.
\end{enumerate}
\end{definition}

\begin{figure}
\begin{center}
\scalebox{0.7}
{
\includegraphics{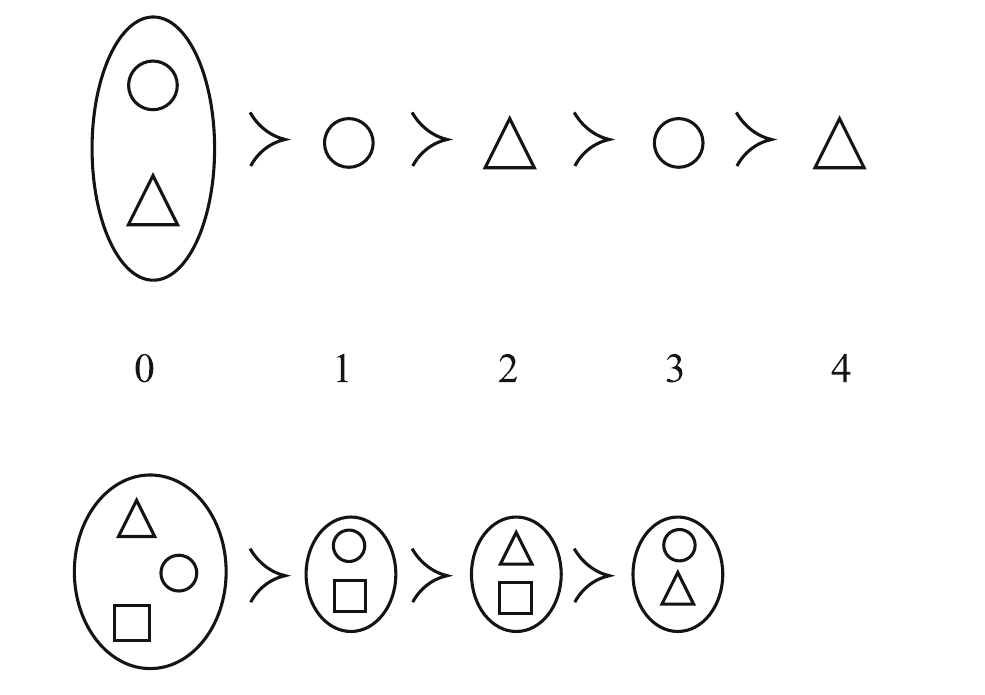}
}
\end{center}
\caption{The models $\mathfrak A(2,2)$ (above) and $\mathfrak A(1,3)$. Points marked by a triangle sa\-tis\-fy $p_1$, circles $p_2$ and squares $p_3$. Ellipses indicate clusters, and the numbers shown represent the $h$-coordinate.}\label{ank}
\end{figure}

\begin{lemm}\label{nkbis} Given natural numbers $N,K,m$, $k\in[1,K]$ and
\[h\in[1,(N-m)K]\]
such that $h\not\equiv k\pmod{K}$,
\[\langle \mathfrak A(N,K),(0,k)\rangle\leftrightarroweq^m_K \langle \mathfrak A(N,K),(h,k)\rangle.\]
\end{lemm}

\proof
Let $\mathfrak A=\mathfrak A(N,K)$, $k\in[1,K]$ and $1\leq h\leq (N-m)K$.

We proceed by induction on $m$. The base case, when $m=0$, is simple, since $(0,k)$ and $(h,k)$ satisfy the same set of propositional variables (namely, $\{p_k\}$).

For the inductive step, we assume that $1\leq h\leq (N-(m+1))K$. Let us first check that $\mathsf{Forth}_\peq$ holds. 

Suppose that $x_1\approx\hdots\approx x_{K-1}\peq (0,k)$. Write $x_i=(\ell,k_i)$ and consider two cases. If $\ell\geq h$, we may set $y_i=x_i$ which clearly satisfy the conditions of $\mathsf{Forth}_\peq$.

Otherwise, there exists a value $k_\ast\in[1,K]$ such that $k_i\not=k_\ast$ for all $i\in[1,K-1]$. Pick $h'\in [(N-(m+1))K,(N-m)K]$ such that $h'\equiv k_\ast\pmod{K}$. For all $i$, $(h',k_i)\in|\mathfrak A|$, and by induction on $m$ we have that\[y_i=(h',k_i)\leftrightarroweq^m_K (0,k_i)\leftrightarroweq^m_K (\ell,k_i),\]
while clearly $y_1\approx\hdots y_{K-1}\peq (h,k),$ hence satisfying $\mathsf{Forth}_\ps$.

For $\mathsf{Back}_\ps$, suppose $y_1\approx\hdots\approx y_{K-1}\peq (h,k)$; then clearly $y_i\peq (0,k)$ for all $i$, and we can set $x_i=y_i$.\endproof

With this, we can show that $\Ps_{i=1}^{k+1}\gamma_i$ cannot generally be defined by formulas of smaller tangled width:

\begin{theorem}\label{morexpressive}
$\mathsf L^{k+1}$ is strictly more expressive than $\mathsf L^k$ for all $k$ over the class of all finite topological models.
\end{theorem}

\proof
Let $\eta^k=\Ps_{i=1}^{k+1} p_i$. Suppose $\varphi\in\mathsf L^k$ has depth $n$ and consider $\mathfrak A=\mathfrak A(n+1,k+1)$.

Then, by Lemma \ref{nkbis}, $(1,2)\leftrightarroweq^{n}_{k+1}(0,2)$, so that by Lemma \ref{basicbis}, $\<\mathfrak A,(1,2)\>\models \varphi$ if and only if $\<\mathfrak A,(0,2)\>\models \varphi$. However, it is easy to check that $\<\mathfrak A,(1,2)\>\not\models \eta^k$ yet $\<\mathfrak A,(0,2)\>\models \eta^k$; hence $\varphi$ cannot be equivalent to $\eta^k$.

Since $\varphi\in\mathsf L^k$ was arbitrary, we conclude that $\eta^k$ is not expressible in $\mathsf L^k$ over the class of finite topological models.
\endproof

\section{Trouble formulas}\label{formulas}

In this section we shall introduce a sequence of formulas $\langle\formsix^k\rangle_{k<\omega}$ with the property that $\vdash^k \formsix^k$. As we shall see later, $\not\vdash^k  \formsix^{k+1}$, thus establishing that $\mathsf{DTL}^{k+1}$ is stronger than $\mathsf{DTL}^k$. The formulas $\formsix^k$ will all be in $\mathsf L^1$.

\begin{definition}
The following abbreviations shall be used throughout the text.
\begin{align*}
\formtwo^k&=\ps p_k\to \displaystyle\bigwedge_{i=1}^k(p_i\to\nx p_{|i+1|_k})\\
\formthree^k_i&=p_i\wedge\hf\formtwo^k\\
\formfour^k&=\nc\displaystyle\bigwedge_{i=1}^k \ps\formthree^k_i\\
\formfive^k&=\displaystyle\Ps_{i\in [1,k]}\formthree^k_i\\
\formsix^k&=\formfour^k\to\hf\ps p_k
\end{align*}
\end{definition}

Before continuing, let us give some intuition for these formulas. The formula $\formtwo^k$ states that $f$ `cycles' the values of $p(x)$; if $p(x)=p_i$, $p(f(x))=p_{i+1}$, unless $i=k$ in which case $p(f(x))=p_1$. The formula $\ps p_k$ is used as a sort of trigger for this cycling behavior; when $\ps p_k$ fails, $p(f(x))$ is unspecified.

$\formthree^k_i$ is used to begin the cycling behavior described by $\formtwo^k$ at $p_i$; it says that, initially, $p_i$ holds, and from then on, $f$ cycles the values of $p(f^n(x))$, provided that $\ps p_k$ holds at each step.

$\formfour^k$ and $\formfive^k$ are similar, but $\formfour^k$ is stronger. As we will mainly be interpreting these formulas over finite Kripke models, let us restrict the discussion to this setting. Here, the meaning of $\formfive^k$ should be familiar; it says there is a cluster where there is a point $x_i$ satisfying each $\formthree^k_i$.

The formula $\formfour^k$, meanwhile, says that each $\formthree^k_i$ is dense near $x$; in particular, each minimal cluster must have one point satisfying each $\formfour^k$. But such a cluster would be a witness to $\formfive^k$.

Thus we have that $\models\formfour^k\to\formfive^k$; but note that the former formula is in $\mathsf L^1$, while the latter is not.

Meanwhile, we should also expect $\models\formfour^k\to\nx\formfour^k$; this is because, if $x_1\approx x_2\approx \hdots\approx x_k$ is a cluster with $x_i$ satisfying $\formthree^k_i$, then clearly each $x_i$ satsifies $\ps p_k$ (since $x_k\peq x_i$), so that $f(x_i)$ satisfies $p_{|i+1|_k}$.

Thus also $f(x_1)\approx f(x_2)\approx\hdots\approx f(x_k)$ is a cluster of points satisfying each $\formthree^k_i$ (although `shifted' one step). It then follows that these points also satisfy $\formfive^k$. By induction, we see that $\models\formfive^k\to\hf\formfive^k$; but this clearly makes $\formsix^k$ true, since $\models\formfive^k\to\ps p_k$.

The reasoning we have just made is easy to formalize in $\mathsf{DTL}^k$, as we show below. Later we shall also see that it is impossible to formalize in $\mathsf{DTL}^{k-1}$.

\begin{prop}\label{isderive}
Given $k<\omega$, $\vdash^k\formsix^k$.
\end{prop}

\proof
Reasoning within $\mathsf{S4}$ one readily sees that, for any $i<k$, $\vdash^0\formfour^k\to\ps(\formthree^k_i\wedge\formfour^k);$
thus we may apply necessitation and $\mathsf{Ind}_\ps$ to derive ${\formfour^k\to\Ps_{i=1}^k\formthree^k_i}$ and obtain
\begin{equation}\label{item1}
\vdash^0\formfour^k\to\formfive^k.
\end{equation}

Further, we note that
\begin{equation}\label{item2}
\vdash^0 \formfive^k\to\ps p_k,
\end{equation}
since this is a consequence of the axiom $\mathsf{Fix}_\ps$.

For any $i\in[1,k]$ we may use $\mathsf{Fix}_\ps$ to see that \[\vdash^0\formfive^k\to\ps(\formthree^k_i\wedge\formfive^k).\]
Using (\ref{item2}), this imples \[\vdash^0\formfive^k\to\ps(\formthree^k_i\wedge\ps p_k\wedge\formfive^k),\]
i.e.
\[\vdash^0\formfive^k\to\ps(p_i\wedge \hf\formtwo^k\wedge\ps p_k\wedge\formfive^k).\]

Now, by $\mathsf{Fix}_\hf$, $\vdash^0\hf\formtwo^k\to(\formtwo^k_i\wedge\nx\hf\formtwo^k)$, whereas
\[\vdash^0 \ps p_k\wedge\formtwo^k\to (p_i\to \nx p_{|i+1|_k}),\]
i.e. $\vdash^0p_i\wedge\ps p_k\wedge\formtwo^k\to \nx p_{|i+1|_k}.$
From this we conclude that
\[\vdash^0\formfive^k\to\ps(\nx p_{|i+1|_k}\wedge \nx\hf\formtwo^k\wedge\formfive^k),\]
and since it holds for all $i\in[1,k]$ we can use $\mathsf{Ind}_\ps$ to obtain
\[
\vdash^0\formfive^k\to\Ps_{i\in[1,k]}(\nx p_{|i+1|_k}\wedge\nx \hf\formtwo^k),
\]
which, rearranging indices and pulling out $\nx$, shows that
\[
\vdash^0\formfive^k\to\Ps_{i\in[1,k]}\nx(p_i\wedge\hf\formtwo^k).\]
Now, we may use $\mathsf{Cont}^k$ to obtain
\[\vdash^k\formfive^k\to\nx\Ps_{i\in[1,k]}(p_i\wedge\hf\formtwo^k);\]
by necessiation and $\mathsf{Ind}_\hf$ this yields \[\vdash^k\formfive^k\to\hf\Ps_{i\in[1,k]}(p_i\wedge\hf\formtwo^k),\]
i.e.
$\vdash^k\formfive^k\to\hf\formfive^k.$

Putting this together with (\ref{item2}) we see that
\begin{equation}\label{nearend}
\vdash^k\formfive^k\to\hf\ps p_k,
\end{equation}
which together with (\ref{item1}) gives us
\[\vdash^k\formfour^k\to\hf\ps p_k,\]
i.e. $\vdash^k\formsix^k$, as claimed.
\endproof

\section{Incompleteness of finite fragments}\label{secinc}

The formula $\formsix^k$ is derivable in $\mathsf{DTL}^k$; let us now see that $\formsix^{k+1}$ is not. To prove this, we shall introduce models $\mathfrak D(N,K)$. They will be composed of two submodels; $\mathfrak C(K)$, defined later, and $\mathfrak B(N,K)$, defined below.

\begin{figure}
\begin{center}
\scalebox{0.7}
{
\includegraphics{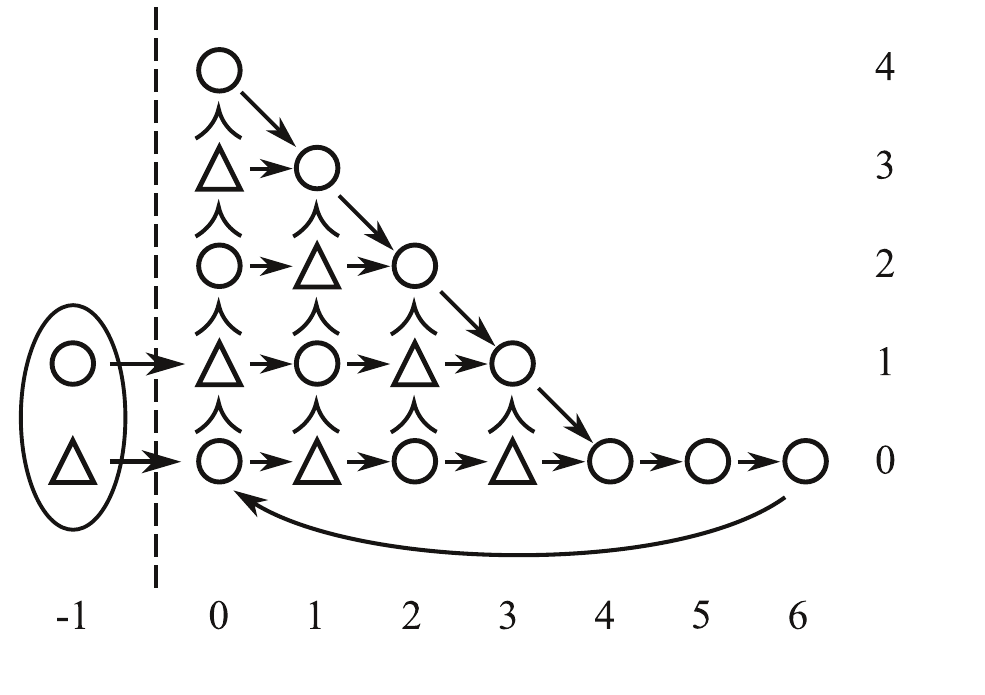}
}
\end{center}
\caption{The model $\mathfrak D(2,2)$, described in Definition \ref{defdnk}. $\mathfrak B(2,2)$, as in Definition \ref{defbnk}, is the submodel on the right-hand side of the dotted line and $\mathfrak C(2)$ is the submodel on its left. Arrows indicate $f_\mathfrak D$, while $\peq$ is the transitive, reflexive closure of the relation represented by $\prec$ together with the ellipse on the left, which represents $\approx$. Points represented by a circle satisfy $p_1$, by a triangle, $p_2$.}\label{fig1}
\end{figure}

The general idea is that the models $\mathfrak D(n+1,k+1)$ will satisfy \[\mathsf{DTL}^k_n\cup \{\neg\formsix^{k+1}\},\]
thus showing that $\not\vdash^k_n\formsix^{k+1}$ for all $n$. From this we may conclude that $\not\vdash^k\formsix^{k+1}$.

Before defining our structures formally, let us give a general idea. Consider the model $\mathfrak D=\mathfrak D(2,2)$ depicted in Figure \ref{fig1}. We will name a point $x$ using triples $(h(x),t(x),k(x))$, where $h(x)$ is the `spatial' (vertical) coordinate, $t(x)$ the `temporal' (horizontal) coordinate and $k(x)$ is the index of $p(x)$, which in this case is $1$ for points represented by a circle and $2$ for triangles. The points on the left of the dotted line will be written $(0,-1,k(x))$.

First, let us observe that $\mathfrak D\models\formtwo^2$, since $f_\mathfrak D$ alternates between circles (which satisfy $p_1$) and triangles (which satisfy $p_2$). The exception for this are the points on the main diagonal $h+t=6$ and on the `tail' $t\geq 4$, but these points do not satisfy $\ps p_2$ and thus they also satisfy $\formtwo^2$. From this, one can easily check that $(0,-1,1)$ satisfies $\neg\formsix^2$.

Meanwhile, the key aspect of the model is that $f_\mathfrak D$ is discontinuous, since $(0,-1,2)\peq (0,-1,1)$ yet
\[f_\mathfrak D(0,-1,2)=(0,0,1)\not\peq(1,0,2)=f_\mathfrak D(0,-1,1).\]

This discontinuity is easily seen to make the following instance of $\mathsf{Cont}^2$ fail on $(0,-1,1):$
\[\ps \nx\{p_1,p_2\}\to\nx\ps\{p_1,p_2\}.\]
However, instances of $\mathsf{Cont}^1$ of small modal depth do hold. Consider, for example,
\[\ps\nx p_1\to\nx\ps p_1.\]
Here we see that $f_\mathfrak D(0,-1,2)$ satisfies $\nx p_1$, so that $(0,-1,1)$ satisfies $\ps\nx p_1$. If $f_\mathfrak D$ were continuous, we would be able to use $f_\mathfrak D(0,-1,2)$ as a witness that $(0,-1,1)$ satisfies $\nx\ps p_1$, but in this case we cannot. However, we do have a different witness, namely $(2,0,1)$. More generally, as we shall see in Lemma \ref{bislemm}, $(2,0,1)\leftrightarroweq^1_\ast (0,0,1)$ so the two satisfy the same formulas of modal depth one.

Thus $\mathfrak D$ satisfies $\mathsf{DTL}^1_1$ as well as $\neg\formsix^2$, from which we conclude that $\not\vdash^1_1\formsix^2$. To see that $\not\vdash^1_n\formsix^2$, we need to consider a larger model, $\mathfrak D(n,2)$, which is built much like $\mathfrak D(2,2)$ but is deeper. By varying $n$, we conclude that $\not\vdash^1\formsix^2$.

Now, let us give the formal definition of $\mathfrak B(N,K)$, which is the submodel of $\mathfrak D(N,K)$ on the right of the dotted lines in Figure \ref{fig1}.

\begin{definition}\label{defbnk}
Given $N,K<\omega$, we define a $K$-simple dynamic model $\mathfrak B=\mathfrak B(N,K)$ by letting
\begin{enumerate}
\item $|\mathfrak B|$ be the set of all triples of natural numbers $(h,t,k)$ such that either
\begin{enumerate}
\item $h+t\leq NK$, $k\in[1, K]$ and $k\not\equiv h+t\pmod K$ or
\item $h=0$, $t\in[NK+1,N(K+1)]$ and $k\not=K$.
\end{enumerate}
\item $(h_1,t_1,k_1)\peq_{\mathfrak B}(h_2,t_2,k_2)$ if and only if $t_1=t_2$ and $h_1\geq h_2$;
\item 
$f_\mathfrak B(h,t,k)=
\begin{cases}
(h,t+1,|k+1|_K)&\text{if $h+t<NK$}\\\\
(h-1,t+1,k)&\text{if $h+t=NK$ and $h>0$}\\\\
(h,t+1,k)&\text{if $t\in[NK+1,N(K+1))$}\\\\
(0,0,|k+1|_{K})&\text{if $t=N(K+1)$}
\end{cases}$
\item $p(h,t,k)=p_k$.
\end{enumerate}
\end{definition}

We will write points as $x=(h(x),k(x),t(x))$. We will also write $s(x)=h(x)+t(x)$.

It will be convenient to describe the $\leftrightarroweq^m_\ast$-equivalence classes over $\mathfrak B(N,K)$. We shall do this using the relations $\sim^m$, defined below.

\begin{definition}
For $m<N$, say $x\sim^m y$ if $p(x)=p(y)$ and one (or more) of the following occurs:
\begin{enumerate}
\item $s(x)=s(y)$,
\item $s(x),s(y)\leq K(N-m)$ or
\item $s(x),s(y)\in[NK,N(K+1)-m].$
\end{enumerate}
\end{definition}

\begin{figure}
\begin{center}
\scalebox{0.7}
{
\includegraphics{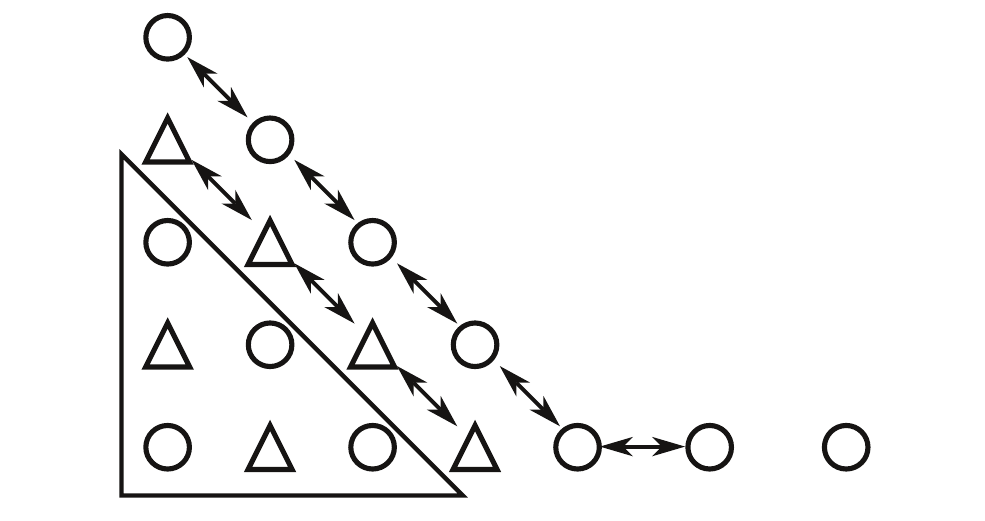}
}
\end{center}
\caption{The relation $\sim^1$ on $\mathfrak B(2,2)$. Points inside the large lower-left triangle are $\sim^1$-related if and only if they satisfy the same propositional variable, which in this figure occurs when they have the same shape. This triangle would grow if we were to consider $\sim^0$, or shrink if we were to consider $\sim^2$.}\label{fig2}
\end{figure}

The models $\mathfrak B(k,n)$ are designed to be very homogeneous, so that different points are hard to distinguished using $\mathsf L^\ast$. The relations $\sim^m$ are representative of this.

First, note that every point is $\sim^m$-similar to another on the $t$-axis; this can be seen in Figure \ref{fig2}, where every point can be `slid' down the diagonals to one on the line $h=0$. Points in the large left-hand triangle may have more than one such representative.

Another useful property is that all points have a very similar orbit; namely, if $x$ is any point and $y$ satisfies $h(y)=0$, then $y$ lies in the orbit of $x$. This is easily seen in Figure \ref{fig1}, where if we follow the $f_\mathfrak D$-arrows starting anywhere, we eventually reach $(0,0,1)$, and from here $f_\mathfrak D$ simply cycles around the $t$-axis indefinitely.

\begin{figure}
\begin{center}
\scalebox{0.7}
{
\includegraphics{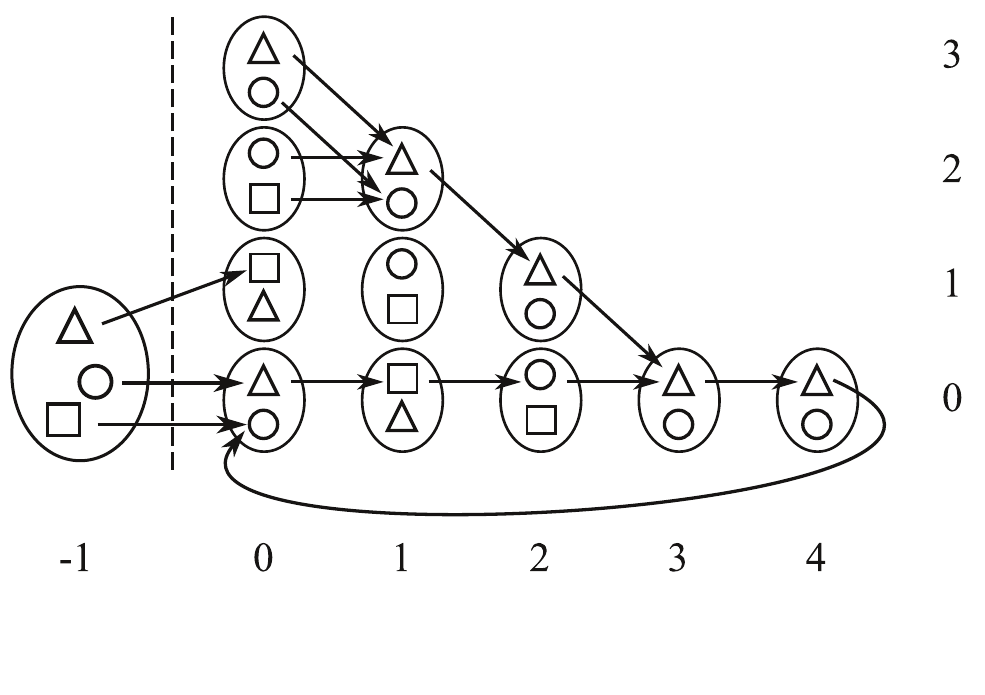}
}
\end{center}
\caption{The model $\mathfrak D(1,3)$. Notation is similar to that for Figure \ref{fig1}, but this time we have only partially indicated the relations; points marked with a square now satisfy $p_3$.}\label{fig3}
\end{figure}

The situation is slightly more involved for larger values of $k$, where we may turn to Figure \ref{fig3}. Consider, for example, the point $(0,0,2)$; in this example, circles have third coordinate $1$, triangles have $2$, squares $3$. Here, notice that after four iterationis of $f_\mathfrak B$ we reach $(0,4,1)$, which then maps to $(0,0,1)$. Afterwards, $f_\mathfrak B$ will cycle through the second row of the $t$-axis, and then return again to $(0,0,2)$.

Let us collect these observations into a lemma:

\begin{lemm}\label{mainaxis}
For every $x\in|\mathfrak B(N,K)|$ and $m<N$,
\begin{enumerate}
\item there is $y\sim^m x$ with $h(y)=0$ and\label{mainaxis1}
\item if $h(y)=0$ there is $n<\omega$ such that $f^n_\mathfrak B(x)=y$.\label{mainaxis2}
\end{enumerate}
\end{lemm}

\proof
The first claim is obvious if we notice that
\[(h,t,k)\sim^m (0,h+t,k).\]

For the second, first we observe that $h(f^{N(K+1)+1}(x))=0$ independently of $x$; then note that $f_\mathfrak B$ is clearly transitive on those elements $z$ with $h(z)=0$, given that
\[f^{N(K+1)+1}(0,0,k)=(0,0,|k+1|_{K}),\]
thus `rotating' $k(z)$.
\endproof

Now, let us see that $\sim^m$ does, indeed, guarantee partial bisimulation.

\begin{prop}\label{bislemm}
If $x\sim^m y$ then $x\leftrightarroweq^m_\ast y$.
\end{prop}

\proof
We work by induction on $m$, considering each clause of a tangled bisimulation. Note that $\sim^m$ preserves atoms, in particular covering the case $m=0$.

Otherwise, suppose $x\sim^{m+1} y$. Clearly we only need to prove the `forth' clauses, since the `back' clauses are symmetric.
\begin{description}
\item[$\mathsf{Forth}_\peq$] We shall only consider the case where $s(x),s(y)< NK$; the other case is similar and easier.

Suppose $x_0\approx x_1\approx\hdots\approx x_{I-1}\peq x$; note that we can assume $I\leq K$, since $\mathfrak B$ has cluster width $K$. Note also that each $x_i$ has $h(x_i)\geq h(x)$ and $t(x_i)=t(x)$.

Consider $h'=h(x_i)+t(y)-t(x)$. If $h'\geq h(y)$, set $h=h'$; otherwise, let $h$ be the least value such that $h\geq h(y)$ and $h+t(y)\equiv h(x_i)+t(x)\pmod{K}$. Then, set $y_i=(h,t(y),k(x_i)).$

First, note that $s(y_i)\equiv s(x_i)\pmod K$, so that all $y_i$ are elements of $|\mathfrak B|$. Now, we further have that $s(y_i)=s(x_i)$ except in the case that $h'<h(y)$, in which it easily follows that $s(x)\not=s(y)$, so $s(x),s(y)<K(N-(m+1))$ and thus $s(x_i),s(y_i)<K(N-m)$.

In either case we use our induction hypothesis to see that $y_i\leftrightarroweq^m x_i$, as claimed.

\item[$\mathsf{Forth}_f$] This follows from observing that the required (in)equalities are preserved by $f_\mathfrak B$ and we skip it.
\item[$\mathsf{Forth}_\hf$] Let $n<\omega$ and consider $z=f^n_\mathfrak B(x)$. Then, by Lemma \ref{mainaxis}.\ref{mainaxis1}, there is $z'\sim^m z$ with $h(z')=0$, while by Lemma \ref{mainaxis},\ref{mainaxis2}, there is $n'$ such that $f^{n'}_\mathfrak B(y)=z'$, as required.
\end{description}
\endproof

Now that we have studied the models $\mathfrak B(N,K)$, let us add the `head' $\mathfrak C(K)$, which is where the trouble really lurks. The resulting model will be called $\mathfrak D(N,K)$, where points in $\mathfrak C(K)$ will map discontinuously onto $\mathfrak B(N,K)$. However, these discontinuities will require large formulas to capture in $\mathsf{L}^{K-1}$, given that $\mathfrak C(K)$ will consist of a cluster with $K$ points.

\begin{definition}\label{defdnk}
We define a model $\mathfrak C=\mathfrak C(K)$ where
\begin{itemize}
\item $|\mathfrak C|=\{0\}\times\{-1\}\times[1,K]$
\item $\peq_\mathfrak C$ is total (i.e., $\mathfrak C$ consists of a single cluster)
\item $p(0,-1,k)=p_k$.
\end{itemize}
We define a model $\mathfrak D=\mathfrak D(N,K)$ based on $|\mathfrak C(K)|\cup|\mathfrak B(N,K)|$ with
\[\peq_\mathfrak D=\peq_{\mathfrak C(K)}\cup\peq_{\mathfrak B(N,K)}\]
and
\[f_\mathfrak D(0,-1,k)=
\begin{cases}
(0,0,|k+1|_K)&\text{if $k\not=K-1$}\\\\
(1,0,K)&\text{if $k=K-1$.}
\end{cases}
\]
\end{definition}

Our strategy now is to show that $\mathfrak D(N+1,K+1)$ is a model of $\mathsf{DTL}^K_N\cup\{\neg\formsix^{K+1}\}$; from this we may conclude that $\not\vdash^K\formsix^{K+1}$, given that $\mathsf{DTL}^K=\bigcup_{n<\omega}\mathsf{DTL}^K_n$.

\begin{lemm}
$\mathfrak D(N+1,K+1)\models \mathsf {DTL}^K_N$.
\end{lemm}

\proof
All the rules of ${\sf DTL}^K_N$ preserve model validity, so it suffices to check that $\mathfrak D(N+1,K+1)$ satisfies all axioms of ${\sf DTL}^K_N$; that is, all permitted substitution instances of axioms of ${\sf DTL}^K$.

Since $\mathfrak D(N+1,K+1)$ is a weak dynamical system, it satisfies every axiom of $\mathsf{DTL}_N^K$ except possibly for instances of $\mathsf{Cont}^K$.

So, let
\[\sigma=\Ps_{i\leq K}\nx\delta_i\to \nx\Ps_{i\leq K}\delta_i\]
be a substitution instance of $\mathsf{Cont}^K$ where each $\delta_i$ has modal depth at most $N$.

Let $x\in|\mathfrak D|$ and assume that
\[x\in\lb\Ps_{i\leq K}\nx\delta_i\rb_\mathfrak D;\]
since $f_\mathfrak D\upharpoonright|\mathfrak B|$ is continuous, we can suppose that $x\in|\mathfrak C|$, for otherwise $x\in\lb\sigma\rb_\mathfrak D$.

Then, given $i<K$ there is $x_i=(0,-1,k_i)\approx x$ such that $x_i\in\lb\nx\delta_i\rb_\mathfrak D$. For at least one value of $k_\ast\in[1, K+1]$ we have that $k_\ast\not=|k_i+1|_K$ for all $i$; we then have that $y_i=(k_\ast,0,|k_i+1|_K)$ is an element of $|\mathfrak B|$ and by Proposition \ref{bislemm}
\[y_i\leftrightarroweq^N_\ast f_\mathfrak D(x_i).\]
Meanwhile, $y_i\approx y_j\peq f_\mathfrak D(x)$ for all $i,j$, so that
\[x\in \lb\nx\Ps_{i\leq K}\delta_i\rb_\mathfrak D,\]
as required.
\endproof

\begin{lemm}
Given $K,N<\omega$ and $k\in[1,k]$,
\[\<\mathfrak D(N,K),(0,-1,k)\>\models \neg\formsix^K.\]
\end{lemm}

\proof
Let $\mathfrak D=\mathfrak D(N,K)$.

First, let us show that every $x\in |\mathfrak D|$ satisfies
\[\formtwo^{K}=\ps p_K\to\bigwedge_{k\leq K}(p_k\to \nx p_{|k+1|_K}).\]

If $s(x)\geq NK$, then $x\not\in \lb\ps p_K\rb_\mathfrak D$ and thus $x\not\in \lb\ps p_K\rb_\mathfrak D$. This shows that $x\in \lb\formtwo^K\rb_\mathfrak D$, as required.

Otherwise, letting $y=f_\mathfrak D(x)$, we note by case-by-case inspection that $k(y)=|k(x)+1|_K$, so that $x$ satisfies $p_{k(x)}\to \nx p_{|k(x)+1|_K}$, whereas for $k\not=k(x)$, $x$ satisfies $p_{k}\to \nx p_{|k+1|_K}$ trivially. Thus $\formtwo^K$ holds everywhere, as claimed.

It follows from this, in particular, that $(0,-1,k)$ satisfies $p_k\wedge \hf\formtwo^K$, i.e. $\formthree^K_k$; this shows that $(0,-1,k)$ satisfies $\bigwedge_{1\leq i\leq K}\ps\formthree^K_i$ and, given that $k\in[1,K]$ was arbitrary,
\[\langle\mathfrak D,(0,-1,k)\rangle\models\nc\bigwedge_{i=1}^k \ps\formthree^K_i=\formfour^K.\]

It remains to show that $(0,-1,k)$ satisfies $\ev\nc\neg p_K$; but this follows from the observation that
\[f^{NK+1}_\mathfrak D(0,-1,k)=(0,NK,k')\not\in\lb\ps p_K\rb_\mathfrak D.\]

We conclude that
\[\langle\mathfrak D,(0,-1,k)\rangle\models\formfour^K\wedge \neg\hf\ps p_K\equiv\neg\formsix^K,\]
as claimed.
\endproof

The following lemma summarizes our results so far:

\begin{lemm}\label{mainlemm}
For all $k<\omega$, the formula $\formsix^{k+1}\in \mathsf L^1$ is derivable in $\mathsf{DTL}^{k+1}$, but not  in $\mathsf{DTL}^{k}$.
\end{lemm}

\proof
By Proposition \ref{isderive}, $\vdash^{k+1}\formsix^{k+1}$; meanwhile, if $\vdash^k\formsix^{k+1}$, we would have that $\vdash^k_n\formsix^{k+1}$ for some $n$.

But this cannot be, since we have seen that
\[\mathfrak D(n+1,k+1)\models\mathsf{DTL}^k_n\cup\{\neg\formsix^{k+1}\},\]
and thus $\not\vdash^k_n\formsix^{k+1}$.
\endproof

With this, we may easily prove our main result.

\begin{theorem}\label{maintheo}
Let $\lambda$ be any language such that ${\sf L}\subseteq \lambda\subseteq {\sf L}^\ast$, and let ${\cal DTL}[\lambda]={\cal DTL}^\ast\cap \lambda$.

Similarly, for $k<\omega$, define ${\sf DTL}^{k}[\lambda]={\sf DTL}^{k}\cap\lambda$.

Then, given any natural number $k$, ${\cal DTL}[\lambda]$ is not finitely axiomatizable\footnote{Observe that ${\sf DTL}^{k+1}$ is finitely axiomatizable over ${\sf DTL}^{k}$, but it does not necessarily follow from this that ${\sf DTL}^{k+1}[\lambda]$ is finitely axiomatizable over ${\sf DTL}^{k}[\lambda]$ for all $\lambda$.} over ${\sf DTL}^{k}[\lambda]$.
\end{theorem}

\proof
Let ${\sf T}$ be any sound, finite extension of  ${\sf DTL}^{k}[\lambda]$, so that without loss of generality we may assume ${\sf T}={\sf DTL}^{k}[\lambda]+\phi$ for some valid formula $\phi$.

Since ${\sf DTL}^\ast$ is complete, we would have that that ${\sf DTL}^\ast\vdash\phi$, and hence, for some value of $K$, ${\sf DTL}^K\vdash\phi$; obviously, we may take $K\geq k$.

But then, we have by Lemma \ref{mainlemm} that ${\sf DTL}^K\not\vdash\formsix^{K+1}$, and hence
\[{\sf DTL}^{k}[\lambda]+\phi\not\vdash\formsix^{K+1}\in{\cal DTL}[\lambda].\] Meanwhile, ${\sf T}$ was arbitrary, so we conclude that ${\cal DTL}[\lambda]$ is not finitely axiomatizable over ${\sf DTL}^{k}[\lambda]$.
\endproof

This result is quite general, so it may be convenient to explicitly mention some special cases. The following corollary states some immediate consequences of Theorem \ref{maintheo};  below, recall that ${\sf KM}\subseteq{\sf DTL}^1[\sf L]$.

\begin{cor}
${\cal DTL}$ and ${\cal DTL}^\ast$ are not finitely axiomatizable. In particular, $\sf KM$ is incomplete for the class of dynamic topological models.
\end{cor}

\section{Concluding remarks}

The axiomatization $\mathsf{DTL}^\ast$ introduced the tangled modality to Dynamic Topological Logic as a sort of scaffolding, to be removed once the appropriate techniques were available. However, I believe the present work to be a convincing argument that indeed it is a central element of the logic; tangled sets affect the behavior of dynamical systems and to be unable to reason about them directly gives a logical formalism an unnecessary handicap.

Of course, none of the results presented here show that a reasonable axiomatization within $\mathsf L^1$ is impossible to find. I am not sure how relevant such an axiomatization would be at this point, but it remains an interesting problem.

Meanwhile, I believe a more fruitful direction is to analyze other logics which are hard to axiomatize because of a similar lack in expressive power. In particular, there are many products of modal logics which have very similar models to those of Dynamic Topological Logic; perhaps they too would benefit from a polyadic variant?

\bibliographystyle{plain}
\bibliography{biblio}

\begin{thebibliography}{10}

\bibitem{arte}
S.N. Artemov, J.M. Davoren, and A.~Nerode.
\newblock Modal logics and topological semantics for hybrid systems.
\newblock {\em Technical Report MSI 97-05}, 1997.

\bibitem{s4real}
Guram Bezhanishvili and Mai Gehrke.
\newblock Completeness of $\mathsf{S4}$ with respect to the real line:
  revisited.
\newblock {\em Annals of Pure and Applied Logic}, 131(1--3):287 -- 301, 2005.

\bibitem{black}
P.~Blackburn, M.~de~Rijke, and Y.~Venema.
\newblock {\em Modal Logic}.
\newblock Cambridge University Press, 2001.

\bibitem{do}
A.~Dawar and M.~Otto.
\newblock Modal characterisation theorems over special classes of frames.
\newblock {\em Annals of Pure and Applied Logic}, 161:1--42, 2009.
\newblock Extended journal version LICS 2005 paper.

\bibitem{dynamictangle}
David~Fern{\'a}ndez Duque.
\newblock Tangled modal logic for topological dynamics.
\newblock {\em Annals of Pure and Applied Logic}, 163(4):467--481, 2012.

\bibitem{me2}
D.~Fern{\'a}ndez-Duque.
\newblock Non-deterministic semantics for dynamic topological logic.
\newblock {\em Annals of Pure and Applied Logic}, 157(2-3):110--121, 2009.
\newblock Kurt G\"odel Centenary Research Prize Fellowships.

\bibitem{me:simulability}
D.~Fern{\'a}ndez-Duque.
\newblock On the modal definability of simulability by finite transitive
  models.
\newblock {\em Studia Logica}, 98:347--373, August 2011.

\bibitem{me:tangle}
D.~Fern{\'a}ndez-Duque.
\newblock Tangled modal logic for spatial reasoning.
\newblock In T.~Walsh, editor, {\em Proceedings of {I}{J}{C}{A}{I}}, pages
  857--862, 2011.

\bibitem{dtlaxiom}
D.~Fern\'{a}ndez-Duque.
\newblock A sound and complete axiomatization for dynamic topological logic.
\newblock {\em Journal of Symbolic Logic}, 2012.
\newblock forthcoming.

\bibitem{konev}
B.~Konev, R.~Kontchakov, F.~Wolter, and M.~Zakharyaschev.
\newblock Dynamic topological logics over spaces with continuous functions.
\newblock In G.~Governatori, I.~Hodkinson, and Y.~Venema, editors, {\em
  Advances in Modal Logic}, volume~6, pages 299--318, London, 2006. College
  Publications.

\bibitem{wolter}
B.~Konev, R.~Kontchakov, F.~Wolter, and M.~Zakharyaschev.
\newblock On dynamic topological and metric logics.
\newblock {\em Studia Logica}, 84:129--160, 2006.

\bibitem{kstrong}
P.~Kremer.
\newblock Strong completeness of $\mathsf{S4}$ wrt the real line.
\newblock 2012.

\bibitem{kmints}
P.~Kremer and G.~Mints.
\newblock Dynamic topological logic.
\newblock {\em Annals of Pure and Applied Logic}, 131:133--158, 2005.

\bibitem{temporal}
O.~Lichtenstein and A.~Pnueli.
\newblock Propositional temporal logics: Decidability and completeness.
\newblock {\em Logic Jounal of the IGPL}, 8(1):55--85, 2000.

\bibitem{mints}
G.~Mints and T.~Zhang.
\newblock Propositional logic of continuous transformations in {C}antor space.
\newblock {\em Archive for Mathematical Logic}, 44:783--799, 2005.

\bibitem{prior}
A.~Prior.
\newblock Time and modality.
\newblock 1957.

\bibitem{tarski}
A.~Tarski.
\newblock Der {A}ussagenkalk\"ul und die {T}opologie.
\newblock {\em Fundamenta Mathematica}, 31:103--134, 1938.

\end{thebibliography}
\end{document}